\documentclass[12pt]{article}
\begin{document}

\title{Computer analysis of the two versions of Byzantine chess}
\author{Anatole Khalfine and Ed Troyan \\
%EndAName
University of Geneva, \textit{\ }\\
24 Rue de General-Dufour, Geneva, GE-1211, \\
Academy for Management of Innovations, \\
16a Novobasmannaya Street, Moscow, 107078\\
e-mail: ed@troyan.me.uk }
\maketitle

\begin{abstract}
Byzantine chess is the variant of chess played on the circular board. In the
Byzantine Empire of 11-15 CE it was known in two versions: the regular and
the symmetric version. The difference between them: in the latter version
the white queen is placed on dark square. However, the computer analysis
reveals the effect of this 'perturbation' as well as the basis of the best
winning strategy in both versions.
\end{abstract}

.\sloppy\tabcolsep=0.15cm

\section{Introduction}

Byzantine chess \cite{BChessPage}, invented about 1000 year ago, is one of
the most interesting variations of the original chess game Shatranj. It was
very popular in Byzantium since 10 CE A.D. (and possible created there).
Princess Anna Comnena\ \cite{Anna} tells that the emperor Alexius Comnenus
played 'Zatrikion' - so Byzantine scholars called this game. Now it is known
under the name of \textit{Byzantine chess}.

Zatrikion or Byzantine chess is the first known attempt to play on the
circular board instead of rectangular. The board is made up of four
concentric rings with 16 squares (spaces) per ring giving a total of 64 -
the same as in the standard 8x8 chessboard. It also contains the same pieces
as its parent game - most of the pieces having almost the same moves. In
other words divide the normal chessboard in two halves and make a closed
round strip \cite{BChessPage}.

\section{Analysis}

It is convenient to imagine it in the form of rectangular diagram \cite%
{Java1}, so that.

\begin{equation}
\begin{tabular}{|l|l|l|l|l|l|l|l|l|l|l|l|l|l|l|l|}
\hline
\phantom{$\bullet$} &  & P & Q & K & P &  &  &  &  & p & k & q & p &  & %
\phantom{$\bullet$} \\ \hline
\phantom{$\bullet$} &  & P & B & B & P &  &  &  &  & p & b & b & p &  & %
\phantom{$\bullet$} \\ \hline
\phantom{$\bullet$} &  & P & N & N & P &  &  &  &  & p & n & n & p &  & %
\phantom{$\bullet$} \\ \hline
\phantom{$\bullet$} & \phantom{$\bullet$} & P & R & R & P & %
\phantom{$\bullet$} & \phantom{$\bullet$} & \phantom{$\bullet$} & %
\phantom{$\bullet$} & p & r & r & p & \phantom{$\bullet$} & %
\phantom{$\bullet$} \\ \hline
\end{tabular}
\label{1}
\end{equation}%
This was obtained from the standard 8x8 board separated into two halves: the
first including the Queen-side and the second including the King-side. Then,
the two halves are glued at their short edges and a circular strip is
obtained: the width of 4 cells and the perimeter of 16 cells. In fact, it is
no more than 4 circular roads of the old Hippodrome of Constantinople,
since, according to the Roman tradition, the races were competed between the
four teams: albati (white), russati (red), prasinati (green) and veneti
(blue).

Now in the Byzantine chess the pieces move as they do in Shatranj - the most
ancient chess variant. Namely, \textsc{kings}, \textsc{rooks}, \textsc{%
knights} obey the standard orthodox FIDE chess rules, but \textsc{bishops}
jump two square diagonally (resemblance with \textsc{knights}), \textsc{%
queens} move one diagonal, \textsc{pawns} have no double first step, nor en
passant, neither they are promoted. \textsc{Pawns} can move clockwise and
counterclockwise so that two pawns of a player going in different directions
may meet on opposing squares, thus blocking each other (however, the late
Byzantines used to treated it as self-annihilation, and the player would
immediately loose both pawns, without counting a move).

The win is achieved by mating the opponent, by stalemating the opponent, or
by `bare \textsc{king}': taking the last piece of the losing side. In the
latter case, the game may reduce to a draw if the losing side, in turn,
manages to capture the last piece of the opponent and, hence, nothing
remaining but the two bare \textsc{kings}. However, mating is not a trivial
task because \textsc{queens} and \textsc{bishops} are very week and their
relative value \cite{Theory} is equivalent to 1.5 pawn (while it is 10 and
3.5 in orthodox FIDE chess \cite{Chessbase}).

Circular chess \cite{Cchess} is a modern game, derived from the original
Byzantine chess. It is played on the circular board with the same notation
and the same initial setup, but the pieces obey the standard FIDE rules
(instead of the Shatranj rules). Indeed, the FIDE\ rules are known to the
widest community and the game of Circular chess has become popular for the
last decade and now it will celebrate the 12 th World Championship. \ This
game involves complicated combinations and analysis although its modern
version was never played in historic Byzantium or Bulgaria of 11-15 CE.

A historic alternative to the regular Byzantine chess \cite{Java1} is its
symmetric version \cite{Java2}%
\begin{equation}
\begin{tabular}{|l|l|l|l|l|l|l|l|l|l|l|l|l|l|l|l|}
\hline
\textsc{\phantom{$\bullet$}} &  & P & \textbf{K} & \textbf{Q} & P &  &  &  & 
& p & k & q & p &  & \phantom{$\bullet$} \\ \hline
\textsc{\phantom{$\bullet$}} &  & P & B & B & P &  &  &  &  & p & b & b & p
&  & \phantom{$\bullet$} \\ \hline
\textsc{\phantom{$\bullet$}} &  & P & N & N & P &  &  &  &  & p & n & n & p
&  & \phantom{$\bullet$} \\ \hline
\textsc{\phantom{$\bullet$}} & \phantom{$\bullet$} & P & R & R & P & %
\phantom{$\bullet$} & \phantom{$\bullet$} & \phantom{$\bullet$} & %
\phantom{$\bullet$} & p & r & r & p & \phantom{$\bullet$} & %
\phantom{$\bullet$} \\ \hline
\end{tabular}%
\label{2}
\end{equation}%
The initial positions of the white \textsc{king} and \textsc{queen} are
replaced so that the white \textsc{king} stands on the light square while
the white \textsc{queen} stands on the dark square: Such replacement results
to the clash between the \textsc{queens}. (Each \textsc{queen} steps one
cell diagonally on the squares of its proper color - as a \textsc{bishop} of
orthodox FIDE chess). Now both \textsc{queens} of symmetric Byzantine chess
pass the same set of squares and, hence, may capture each other. \ 

This symmetric version of Byzantine chess should lead to the change in the
game tactics and strategy. The Queens can encounter (in orthodox FIDE chess\ 
\cite{Euwe}. the endings with similar Bishops are rarely resulting to draw).
Will it imply complication and sharpness of symmetric version? It can be
revealed only in the direct play.

The online Java applets allows to perform the instant analysis. We have
launched one hundred Computer-vs-Computer games for each variant.

In the regular version \cite{Java1} White won \textbf{17} games, lost\textbf{%
\ 13} games, draw was recorded in \textbf{70} games.

In the symmetric version \cite{Java2} the competition yielded the result: 
\textbf{+26=61-13}.

As an alternative method we used the PC program 'Byzantine Chess' based on
Zillions of Games platform \cite{Zill1, Zill2}. We launched 20 games in
Computer-vs-Computer competition with 20 sec time control per move.

In the regular version White won \textbf{5} times, lost\textbf{\ 4} times,
draw was recorded in \textbf{11} games.

In the symmetric version the competition yielded the result: \textbf{+5=12-3}%
.

In fact, to demonstrate the smoothness of the Byzantine chess with respect
to the Circular chess, we launched the latter \cite{Java3} in a 100-game
Computer-vs-Computer match, that scored as \textbf{+42=34-24 } (of course
White won).

\section{Conclusion}

According to computer analysis, the symmetric version of Byzantine chess 
\cite{Java2} reveals a bit less 'peaceful policy' with respect to the
regular version \cite{Java1}. We could recognize the similarity with the 
\textsc{bishop-vs-bishop }endgames of orthodox FIDE chess \cite{Euwe}.
Namely, the endings with dark-vs-light pieces usually admit 'peaceful'
coexistence, but dark-vs-dark and ligth-vs-light pieces cannot live without
a 'war'. In fact, the latter case corresponds to the symmetric version of
Byzantine chess.

We have also compared the Byzantine chess and Circular chess (played on the
same round board but according to the FIDE rules), and discovered that in
the latter variant draws occur approximately two times less frequent (1/3
and 2/3 of all games, respectively). Large number of draws is not strange
for the Byzantine chess. There is no \textsc{pawn} promotion and many
endgames are played upto the stage of a bare \textsc{king} (the condition of
loss). The game tends to draw when the major pieces (\textsc{rook} and often 
\textsc{knight}) are removed off the board - the remaining \textsc{pawns}, 
\textsc{bishops} and \textsc{queens} bring no tension in spite of possible
extra piece belonging to the opponent. Moreover, \textsc{Knight-vs-Bishop}
and \textsc{Knight-vs-Queen} endgames do not promise any chance to win, for 
\textsc{bishop} and \textsc{queen} can provide perpetual survival against
the stronger piece. However, a small deviation from equilibrium in \textsc{%
Pawn-vs-Pawn} endgame \cite{Queen} often results to drastic circumstances.
As for mating, it is possible only in the opening or, at least, in the
middlegame, especially, if there is some difference between the skill of the
players. Computer-vs-Computer program has considered this situation only 3
times in 200 games, although a serious player will mate any existing online-
or, offline-program: the Elo rating of the Byzantine chess engines is still
lower the \ FIDE master (FM) level.

Much more strange that White has advantage over Black. It is absolutely
evident in the symmetric version of Byzantine Chess. Have we discovered the
specific effect resulted from the \textsc{queens} belonging to the same
color (e.g. both move on the dark squares)? We cannot explain this strange
fact, perhaps it is due to the peculiarity of the chess engine, while human
player may never encounter this effect.

Anyhow, the best game strategy is now outlined: achieve advantage in the
middlegame, exchange all minor pieces (\textsc{bishop}, \textsc{queen}) and
reduce to \textsc{pawns} and major pieces (\textsc{rook}, \textsc{knight}) \
in the endgame. Otherwise, there would be no use of one or two extra minor
pieces. On the other hand, it is possible to avoid the evident loss by
exchange of all major pieces and \textsc{pawns}, thus, reducing to the only
minor pieces (but not \textsc{pawns} instead!).\textbf{: }

\section{Perspective}

If we compare Byzantine chess with Circular chess, \ - that is a
Byzantine-FIDE hybrid of usual chess played on the circular board - we see
that the Byzantine chess is a quiet game with plain development and absence
of complicated combinations (perhaps pertaining to the style of live and
mind of the late Byzantine scholars and noblemen) - no wonder that stalemate
and bare \textsc{king} are defined as loss. Mating the opponent's \textsc{%
king} is so brilliant event that we believe it should be rewarded 1.5 points
(instead of standard 1 point).

Another note concerning the rule of draw: many endings are deprived of real
possibility to win and it is reasonable to stop the game if there is no
capture, say, within 20 moves. The evident draw endgames (\textsc{%
Bishop-vs-Queen}, \textsc{2Bishops-vs-Bishop\&Queen} etc.) should be also
stopped without further play.

In spite of the prosperous development of Circular chess \cite{Cchess} in
the last decade, the native original Byzantine chess is played sporadically.
The organization of regular tournaments is still under discussion, and it is
important to regulate the rules. What we cannot know is the exact scores of
Byzantine chessmasters but we can discover and restore the theory of game.


\begin{thebibliography}{99}
\bibitem{BChessPage} H.L. Bodlaender Byzantine chess. A variant of Shatranj,
played on a round board (1997). \
http://www.chessvariants.org/historic.dir/byzantine.html

\bibitem{Anna} The Alexiad of Anna Comnena. Book XIV. (Penguin Books, 1969).

\bibitem{Java1} E. Friedlander \ Byzantine Chess Java-Applet. \
http://play.chessvariants.org/erf/Byzantin.html

\bibitem{Theory} A. Khalfine and E. Saperow. \ Towards the theory of
Byzantine Chess. (2002). \
http://www.chessvariants.org/historic.dir/byztheory.html

\bibitem{Chessbase} ChessBase Platform. \ http://www.chessbase.com

\bibitem{Cchess} Circular Chess Society. \ http://www.circularchess.org

\bibitem{Java2} E. Friedlander \ Byzantine Chess Java-Applet, Symmetric
version.. http://play.chessvariants.org/erf/Byzanti2.html

\bibitem{Euwe} M. Euwe and H. Kramer. \ \ The Middlegame - Book I : Static
Features (Hays Pubs 1994)

\bibitem{Zill1} R. Price. Byzantine Chess (2001). PC Program on the
Zillions-of-games platform.

http://www.chessvariants.org/index/zillions.php?itemid=zByzantineChess

\bibitem{Zill2} Anatoly Khalfine and Ernst Saperow. Toroidal Byzantine Chess
Rectangular (2002). PC Program on the Zillions-of-games platform.
http://www.zillionsofgames.com/cgi-bin/zilligames/submissions.cgi/22549?do=show;id=926

\bibitem{Java3} E. Friedlander \ Circular Chess Chess Java-Applet. \
http://play.chessvariants.org/erf/Circular.html

\bibitem{Queen} E. Saperow. \ Byzantine Chess endgames. Queen vs Pawn.
(2003)
http://www.chessvariants.org/historic.dir/byztheory.files/byzendqp.html
\end{thebibliography}
\end{document}